%% file: localcuts.tex
\RequirePackage{fix-cm}
\documentclass[smallextended]{svjour3}       
\smartqed  
\usepackage{graphicx}

\input{packages.tex}
\input{commands.tex}

\begin{document}

\title{Learning to Use Local Cuts}

\author{Timo Berthold         \and
        Matteo Francobaldi \and
	Gregor Hendel
}

\institute{T.\ Berthold, G.\ Hendel \at
Fair Isaac Germany GmbH \\
Takustr. 7, 14195 Berlin, Germany \and
M.\ Francobaldi \at 
Department of Computer Science and Engineering, University of Bologna \\
Viale Risorgimento 2, 40136 Bologna, Italy
}
\date{Received: date / Accepted: date}

\maketitle

\begin{abstract}
An essential component in modern solvers for mixed-integer (linear) programs (MIPs) is the separation of additional inequalities (cutting planes) to tighten the linear programming relaxation. Various algorithmic decisions are necessary when integrating cutting plane methods into a branch-and-bound (B\&B) solver as there is always the trade-off between the efficiency of the cuts and their costs, given that they tend to slow down the solution time of the relaxation. One of the most crucial questions is: Should cuts only be generated globally at the root or also locally at nodes of the tree?
We address this question by a machine learning approach for which we train a regression forest to predict the speed-up (or slow-down) provided by using local cuts.
We demonstrate with an open implementation that this helps to improve the performance of the FICO Xpress MIP solver on a public test set of general MIP instances. We further report on the impact of a practical implementation inside Xpress on a large, diverse set of real-world industry MIPs.
\keywords{Mixed Integer Programming \and Cutting Planes \and Machine Learning}
\subclass{90C11 \and 90C05}
\end{abstract}

\section{Introduction}
\label{1}
\emph{Cutting planes} play an indispensable role in solving mathematical optimization problems, and they are at the heart of all competitive solvers for mixed-integer programming today.
Consequently, there is a plethora of scientific papers on cutting planes algorithms, how to improve them, which theoretical characteristics they have, how to improve the numerical stability of cuts, which cuts to select from a given set and so on.
We refer to the excellent surveys~\cite{Marchand2002,Cornuejols2008} and the references therein for an overview of computationally efficient cutting plane methods for general mixed-integer programming.
However, there is surprisingly little work on whether for a given MIP instance cuts should be continuously generated throughout the tree search or whether the tree search would be more efficiently handled by a pure B\&B approach.
The present paper will address  this question and introduce a Machine Learning approach to identify MIP instances that benefit from not applying cutting planes locally inside the tree.

\section{Problem discussion: To cut or not to cut}
\label{2} 
Cutting planes can be integrated into the B\&B scheme in many different ways, each one determined by a combination of several algorithmic decisions~\cite{bcopt,Martin2001,FUGENSCHUH2005,Achterberg2007}.
In this paper, we classify a cutting plane either as a \emph{global cut} or as a \emph{local cut}, according to whether it is generated at the root node or at an internal node of the B\&C search tree.

As for any other routine involved in the solving process of a \MIP, designing an efficient cutting strategy is a matter of trade-offs.
On one side, the use of cuts improves the dual bound, hence increases the chances of pruning, and is consequently expected to reduce the number of explored nodes.
On the other side, cut generation does not come for free; in particular, their inclusion enlarges the size of the LP relaxation, thereby slowing down the linear solver.
Moreover, the use of local cuts, whose validity can not be guaranteed globally, prevents the use of conflict analysis~\cite{ACHTERBERG20074,WitzigBertholdHeinz19}. This might lead to a performance degradation on instances for which conflict analysis is beneficial.

In the present work, we address the question of whether to cut at internal nodes of the B\&B tree, or to limit all cutting activity exclusively to the root node; equivalently, whether to use local cuts during the solve or not.
We denote these two alternatives by $\C$ and $\NC$, respectively.
To the best of the authors' knowledge, up until the time of writing this paper, no efficient criterion has been discussed in the literature to provide an answer to this question for general \MIP.
According to a computational study conducted on our dataset with \solver{FICO Xpress}, the $\C$ method is, on average, 27\% faster than $\NC$.
However, on 22\% of the test set, $\NC$ was significantly faster than $\C$, against the general trend.
If we had a perfect oracle to decide, for any given input instance, whether $\NC$ works better than $\C$, then we could speed up the average runtime of the solver by 11\%.
This statistic indicates the potential for a method that could discriminate between instances that solve better with one or the other approach.
At the same time, it suggests an upper limit on the improvement that we can expect by such a method.
The goal of this paper is to train a machine learning model to predict whether the $\NC$ or the $\C$ approach will work better for a given instance. 

\section{Methodological Approach}
\label{3}
We developed our ML approach by using the benchmark set \miplib{2017}~\cite{MIPLIB2017} as problem set, and \solver{FICO Xpress 8.9} as MIP solver. Note that, however, our methodology can be implemented for any \MIP\ solver $\mathscr{S}$ (supporting both configurations $\C$ and $\NC$) and problem set $\mathcal{P} \subseteq \mathscr{P}$, where $\mathscr{P}$ denotes the space of all MIP problems.

\paragraph{Feature Design.} We represent a problem $p \in \mathscr{P}$ as a vector of $32$ features, which condense the relevant pieces of information leading to an algorithmic discrimination between $\C$ and $\NC$.
Table \ref{Features_Tab} provides a comprehensive description of the resulting \emph{feature space}, denoted by $\mathcal{F}_{\mathscr{S}} \subset \mathbb{R}^{32}$, together with the definition of each individual feature. 
We denote by $A^\prime, b^\prime, c^\prime$ the set of absolute values of the non-zero entries of the problem's matrix $A$, vectors $b$ and $c$, respectively:
$A^\prime := \{\lvert A_{i,j}\rvert; (i,j) \in [m] \times [n], A_{i,j} \neq 0 \}$,
$b^\prime := \{\lvert b_{i} \rvert; i \in [m], b_{i} \neq 0 \}$, and
$c^\prime := \{\lvert c_{j} \rvert; j \in [n],  c_{j} \neq 0 \}.$
The number of remaining rows and columns after presolving, instead, are denoted by $\tilde{m}$ and $\tilde{n}$, respectively.
Finally, we consider the following objective values at the root node: the objective value of the initial LP optimum (initial bound), of the LP optimum at the end of the root node cutting loop (dual bound), and of the incumbent at the end of the root (primal bound), if already available.
We denote them as $\initbound$, $\dualbound$, and $\primalbound$, respectively.

\begin{table}[h!]
        \centering
        \input{Features_Tab_twocol}
        \caption{Our feature space $\mathcal{F}_{\mathscr{S}}$. Note that, for some of these features (i.e., original and presolved rows and columns, as well as order of magnitude of problem data), we use the transformation $\ln()$.
        This transformation, although irrelevant for our $\M$, was indeed convenient for different models that we used in preliminary experiments.}
        \label{Features_Tab}
\end{table}

Features can be classified as either \emph{static} or  \emph{dynamic}. The former represent the mathematical formulation and combinatorial structure of the problem, hence they are completely solver-independent; the latter, instead, capture the algorithmic behavior and development of the problem, hence they are solver-dependent.
Moreover, we divide the static features into three groups, that is, \emph{Matrix}, \emph{Variables} and \emph{Constraints}, corresponding to different components of the problem formulation.
The first group provides information about the size and the sparsity of the problem matrix, as well as the presence of eventual symmetries; the second and the third describe instead the variable and constraint composition of the problem~\cite{MiplibWebsite}, respectively.
The dynamic features are split into $3$ groups as well, that is, \emph{Presolving}, \emph{Scaling} and \emph{Global Cutting}, corresponding to different stages of the solving process, in which the solution process is not yet affected by our decision.
The first group describes the order of magnitude of the problem data after the scaling process, the second measures the effectiveness of the presolver, while the latter quantifies the impact of the global cutting loop. 

\begin{figure}[h!]
\centering
\includegraphics[scale = 0.35]{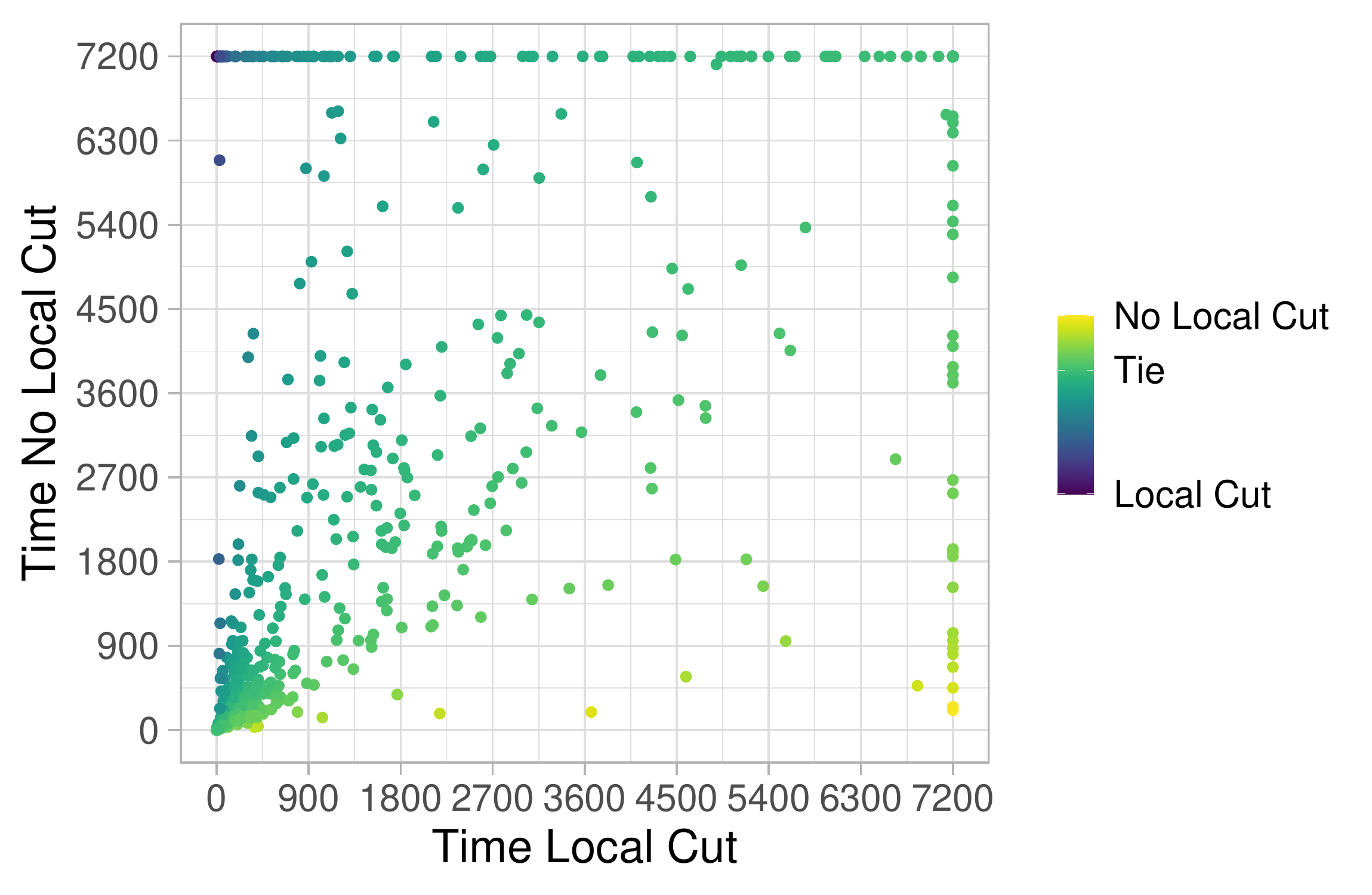}
\caption{The color varies continuously, from the upper-left corner to the lower-right one, between purple, representing those instances on which $\C$ and $\NC$ reach their minimum and maximum runtime, respectively, and yellow, corresponding to the opposite situation. Along the diagonal $y = x$, the color results from the interpolation of the two colors at the extremes, indicating that, on these points, the two methods perform similarly.}
\label{speedup_point}
\end{figure}

\paragraph{Label Definition.} The seemingly natural approach to solve our problem would be to formulate it as a \emph{classification} task, given that the problem itself consists, essentially, in taking a binary decision.
There are, however, two reasons which are in favor of taking a \emph{regression} approach. Firstly, our ultimate goal is to improve the average runtime of the solver, which is a metric that is numerical and not categorical.
Secondly, our focus is on getting the prediction right for those instances on which $\C$ and $\NC$ significantly diverge from each other. For instances where $\NC$ and $\C$ show a similar behavior, the classification will have a negligible impact on the solver performance. 
We define the label of a problem $p \in \mathscr{P}$ as the \emph{speedup factor} between $\mathit{Time}_{\C}(p)$ and $\mathit{Time}_{\NC}(p)$, rescaled by means of the $log_2$ function:
\[
y_{\mathscr{S}}^p = log_{2} \left(\frac{\mathit{Time}_{\C}(p) + 1}{\mathit{Time}_{\NC}(p) + 1}\right) \in \mathbb{R},
\]
where $\mathit{Time}_{\mathcal{\phi}}(p)$ denotes the running time of $\mathscr{S}$ while solving $p$ with configuration $\phi \in \{\C, \NC\}$. The runtimes of the two methods are both augmented by 1, to mitigate the impact of very small numbers, and to prevent the division by zero.
The scatter plot in Figure \ref{speedup_point} displays the instance-by-instance comparison between the runtimes of the two methods over $\miplib{2017}^*$ (discussed in the following paragraph), computed with Xpress and colored continuously according to the speedup value. Now, we can provide a rigorous formulation of our regression problem: our task is to produce a model $M_{\mathscr{S}} \colon \mathcal{F}_{\mathscr{S}} \longrightarrow \mathbb{R}$ that, for each input vector $x_{\mathscr{S}}^p$, approximates the speedup factor $y_{\mathscr{S}}^p$ between the two runs, $\C$ and $\NC$, of the solver $\mathscr{S}$ over $p \in \mathscr{P}$, that is, $M(x_{\mathscr{S}}^p) \approx y_{\mathscr{S}}^p$. In Section~\ref{4}, we will confirm our expectations on the superiority of this regression formulation over a classification one.

\paragraph{Data Collection.} For the given solver $\mathscr{S}$ and problem set $\mathcal{P}$, we collect a ground dataset for our ML approach as follows. Firstly, we apply six random permutations to the elements of $\mathcal{P}$, each one associated to a seed $s = 0, \dots, 5$ (with $0$ identifying the identity permutation), in order to enlarge and diversify our ground problem set. We denote the expanded set by $\mathcal{P}^*$, and each of its elements by $p_s$, referred to as an \emph{instance} of the problem set. Note that the six instances produced from $p$, although mathematically equivalent, can have a very different computational behavior~\cite{LodiTramontani2014}. Hence, from each instance $p_s$, we collect the solver-independent information that we use to compute the static features.

Then, we run the solver $\mathscr{S}$ over $p_s$ twice, once with configuration $\C$, and once with $\NC$. From these runs, we retain the solver-dependent data that we use to compute both the dynamic features and the labels. Now, from the raw information collected, we construct the dataset of features-label observations $\mathcal{D}  = \mathcal{D} \left(\mathscr{S}, \mathcal{P}\right)  = \{(x_\mathscr{S}^{p_s}, y_\mathscr{S}^{p_s}): \,\,\, p_s \in \mathcal{P}^*\} \subseteq \mathcal{F}_{\mathscr{S}} \times \mathbb{R}$, that we use for our learning experiments. In particular, we split this set into a \emph{training set}, $\mathcal{D}_{train}$,  and a \emph{test set}, $\mathcal{D}_{test}$, as follows: we use the permuted problems of $\mathcal{P}$ ($s = 1, \dots, 5$) for the former, and the original ones ($s = 0$) for the latter, corresponding to roughly a training-test split of $83\%$ to $17\%$ of the entire dataset.
Using permutations for a training-test split has first been suggested in~\cite{berthold2021learning}.

\paragraph{Training \& Testing.} Among the different ML models that we considered for our learning task, the best results were provided by a \emph{random forest} (\M) ~\cite{Ho1995}, an ML model consisting of a multitude of \emph{regression trees}~\cite{MorganSonquist,MessengerMandell}, trained by means of ensemble techniques able to improve robustness and prevent overfitting. An overview of the preliminary experiments that we conducted with other ML models can be found in the master's thesis of the second author~\cite{Francobaldi2022}.

We train \M by using the \texttt{caret}~\cite{caret} package of  \texttt{R}~\cite{RLanguage}, with \texttt{method} set to \texttt{"rf"} in the \texttt{"train"} function. Moreover, in order to select a suitable number of regression trees to employ within the ensemble, we use the  \texttt{trainControl} tool with \texttt{method} set to \texttt{"cv"}  and \texttt{number} to \texttt{5}, which performs a $5$-\emph{fold cross-validation} with randomly generated folds~\cite[Section~11.2.4]{UndMachineLearning}. The produced model is a random forest with $500$ regression trees.

We assess the quality of our trained model when deployed in a MIP solver by three different criteria: the running time (\textit{Time}), the primal-dual integral (\textit{PDI})~\cite{BERTHOLD2013611} and the number of explored nodes (\textit{Nodes}). For each of these metrics $\Met$, we compute the performance of \M, on a certain problem $p$, as
\[
\Met_\M(p) =
\begin{cases}
\Met_\C(p), & \text{if} \,\,\, \M(x^p_\mathscr{S}) \leq \tau \\
\Met_\NC(p), & \text{if} \,\,\, \M(x^p_\mathscr{S}) > \tau,
\end{cases}
\]
where $\Met_\phi(p)$ is the performance of our solver $\mathscr{S}$ when solving $p$ with strategy $\phi$. Here, $\tau \in \mathbb{R}$ is a \emph{threshold} that we use as a switch between $\C$ and $\NC$ to convert our regressor into a binary classifier. As experimentally derived, we set the value of $\tau$ to $0$, which corresponds to both methods being predicted to take the same run time.

We aggregate the performance of a strategy $\phi$, over a set of problems $\mathcal{T} \subset \mathscr{P}$, by means of the \emph{shifted geometric mean},
\[
\mathit{Shm}_{\Met}(\phi) = \left(\prod\limits_{p \in \mathcal{T}} \left(\Met_\phi(p) + S\right)\right)^{\frac{1}{\mid \mathcal{T}\mid}} - S,
\]
with \emph{shift} $S$ set to 10 for \textit{Time} and \textit{PDI}, and to 1000 for \textit{Nodes}.

\section{Computational Experiments}
\label{4}

Our computational study consists of two parts. The first experiment was conducted on the dataset $\mathcal{D} = \mathcal{D}\left(\text{Xpress},\miplib{2017}\right)$, collected by executing our data collection process (Section \ref{3}) with version 8.9 of the FICO Xpress MIP solver  and the benchmark set of \miplib{2017} as instance set. We chose a time limit of $T = 7200s$.
As competitors of \M, we consider the following strategies: \always, which always chooses $\C$, \never, which always chooses $\NC$, and \oracle, which always takes the optimal choice, each time according to the particular metric used, between $\C$ and $\NC$, hence representing the best possible performance that we can achieve with our selection method. The main goal for this first experiment was to analyze the potential of the approach and to evaluate the impact of different decisions within our framework.

For the second experiment, we implemented a learned model into FICO Xpress and tested it on an internal data set, which primarily consists of client instances. This implementation is used by default in the current release of Xpress. The main goal of this experiment is to demonstrate the impact of a practical implementation of our approach on a diverse set of real-world instances.

For the first experiment, we excluded from $\mathcal{D}$ all the instances that are not suitable for our study: the ones for which any of the strategies ran out of memory, the ones solved already at the root, and a few cases where one of the strategies ran into numerical issues. This cleaning process reduced the size of $\mathcal{D}$ from 1440 instances to 1155 instances on which the model was trained.


We compare our learned model with its competitors, on both $\mathcal{D}_{train}$ and $\mathcal{D}_{test}$. The results of the evaluation are reported in Table \ref{Regression_tab}, where \textbf{Imp.} refers to the improvement provided to the solver by our model $\M$, when compared against the better of our two competitors \always\ and \never, \textbf{Pot.} denotes the potential improvement, i.e., the performance gap between this competitor and \oracle.
\textbf{Achiev.} reports the percentage of this potential improvement that $\M$ is able to achieve. An immediate observation that we can make is that, between \always\ and \never, the former is certainly the more performant one, given that, on both sets, it significantly outperforms \never in terms of all considered metrics. This confirms our claim from Section \ref{2}: if \always\ and \never\ were the only available strategies to take the $\C$/$\NC$ decision, then we would choose the former.  The results scored by our learned strategy \M, however, shows that a smarter strategy than \always\ is achievable.

\begin{table}[h!]
\input{Regression_tab.tex}
\caption{The \textit{Shm} performance of the competing strategies, in terms of different metrics.}
\label{Regression_tab}
\end{table}

When compared against \always\ in terms of \textit{Time}, \M\ provides a speedup of roughly 8.5\% on the train set and of 3.3\% on the test set, hence achieving about $\nicefrac{3}{4}$ and $\nicefrac{1}{3}$ of the \emph{potential} improvement on the two sets, respectively. Similar considerations can be made by looking at the \textit{PDI}, even though, in this case, the contribution provided to the solver is less pronounced than the one observed in \textit{Time}.
In contrast, the results obtained in terms of \textit{Nodes} seem to contradict the ones observed in the two metrics previously discussed. Indeed, with an improvement of less than 1.5\% and an achievement of less than 15\%, the \textit{Nodes} performance of our solver is, on $\mathcal{D}_{train}$, almost unaffected by the use of our model, while it is even weakened  by it on $\mathcal{D}_{test}$, where we can observe a 7\% degradation and a 70\% increment in the average node consumption.

This tendency, however, is explainable by the fact that the use of cutting planes, as discussed in Section \ref{2}, has the main advantage of reducing the number of explored nodes. This is why, when measured in terms of \textit{Nodes}, the \always\  strategy tends to show the smallest number of nodes. To provide a more detailed evaluation of the contribution given by our approach, we compare our model, \M, against our main competitor, \always, on the entire dataset $\mathcal{D}$ and on different \emph{brackets} of this set, as well as on the subset of the \emph{affected} instances.

The results of this comparison are reported in Table \ref{brackets}.
Here, the bracket $[t_1, t_2]$ is the set of instances of $\mathcal{D}$ satisfying the following: 1) they are solved (within the time limit $T$) by at least one of the two strategies, 2) the runtime of the slower strategy is between $t_1$ and $t_2$. We fix  the right-hand side of each bracket at the time limit $T$, while we increase the left-hand one progressively, so to define a hierarchy of subsets of $\mathcal{D}$ of increasing difficulty. The affected instances, instead, are the ones that are solved by at least one of the two competitors, and for which these competitors make opposite choices, i.e., the ones on which \M decides to deactivate local cutting, as opposed to the default strategy.

\begin{table}[h!]
\centering
\input{Bracket_tab.tex}

\caption{Comparison between \M\ and \always, on different brackets of $\mathcal{D}$.}
\label{brackets}
\end{table}

We observe that the learned model \M\ is able to solve $18$ instances more than \always, and to improve the average runtime by $7.7 \%$. This increases to more than $11\%$ on the instances of the bracket $[10, 7200]$, and keeps increasing together with the hardness of the evaluation set, until reaching an encouraging value of $18.7\%$ on the most difficult among our problems. Finally, by restricting the test bed to the affected instances, we can observe a very promising improvement of more than 26\%, cf.~\cite{Francobaldi2022}.
In an offline experiment, we also tested a classification version of our random forest training. While it also improved the performance, the benefit was clearly smaller than for the regression forest.


\begin{table}[h!]
\centering
\input{Bracket_tab_Xpress.tex}
\caption{Comparison between a version of Xpress with and without a \M\ model, on different brackets of an internal data set of FICO.}
\label{brackets_xpress}
\end{table}

Finally, the results of our initial study were used to implement a new decision module into the FICO Xpress MIP solver~\cite{xpress} for our second experiment.
It was trained and tested on FICO internal data sets which mainly consist of real-world instances provided by FICO customers.
Naturally, this led to a model that is different in detail, but the applied training procedure was exactly the same.
We restrained ourselves to using at most seven features and 50 trees for the final model for three reasons: to save memory, to speed up the \M\ evaluation, and to avoid overfitting.
Both restrictions performed only slightly worse than a variant that used 500 trees and the full feature set.
Given the known average superiority of $\C$ over $\NC$, as an additional safeguard against performance losses, we skewed our decision towards the former approach as follows: in order to deactivate local cuts, we required not only the overall model prediction to be in favor of $\NC$, but also a majority of 70\% of the individual trees.
In other words, we kept generating local cuts even though the random forest decided to deactivate them, unless more than 70\% of the ensemble agreed with this decision.


The results of running a pre-release version of FICO Xpress~8.13 can be seen in Table~\ref{brackets_xpress}.
The presented data set of 5530 test instances constitutes the standard performance evaluation set of Xpress.
Note that it is significantly larger than the \miplib{2017} set we used for our initial experiment, and that this also includes instances which are solved at the root node.
Hence, we expect the results to be less pronounced.
Note further, that this is not the set that was used for training (and testing), but the result of a consequent verification.

We observe an overall improvement of 1.44\% when using a random forest to decide on a local cutting strategy, again with increasing impact as the models get harder.
For instances which take at least 1000 seconds to solve, we see an improvement of 5.25\%.
Additionally, there are nine more instances solved when using an automatic decision on local cuts.
Again, the effect is much more pronounced when we restrict our attention to the affected models, hence those where we deactivate cutting now, but haven't done before. Here, we observe a speedup of almost 18\%, which is comparable to the results of our initial experiment on \miplib{2017}.
Moreover, we used a Wilcoxon signed rank test~\cite{Wilcoxon1945} to evaluate the significance of the improvement.
It rejected the null hypothesis with a p-value of less than 0.001, which can be interpreted as our observations indicating an actual speedup with more than 99.9\% confidence.

\section{Conclusion \& Outlook}

We introduced the first, to the best of our knowledge, methodology to predict for a given general MIP instance whether it would be better solved by a cut-and-branch or a branch-and-cut approach.
Our methodology showed very promising results both for the \miplib{2017} benchmark set and for an internal performance evaluation set of FICO.
Consequently, it has been implemented by the FICO Xpress developer team and is used by default as of the recent 8.13 release of Xpress.

As possible outlook of the present work, we suggest to learn the maximum depth of the B\&B tree in which the cutting procedure should be stopped, rather than choosing only whether to deactivate the procedure after the root node or not.
This seems particularly promising when connected with the concept of lifting local conflicts as suggested in~\cite{WitzigBerthold2021}. 


\bibliographystyle{spmpsci}      
\bibliography{sn-bibliography}


\end{document}

%% file: packages.tex
\usepackage{nicefrac}
\usepackage[scr=boondoxo]{mathalfa}
\usepackage{adjustbox}
\usepackage{colortbl}
\usepackage{hhline}
\usepackage{verbatim}
\usepackage{amsmath}
\usepackage{amssymb}
\usepackage{array}
\usepackage{xspace}
\usepackage{multirow}
\usepackage{doi}

%% file: commands.tex
\newcommand{\MIP}{\text{MIP}\xspace}

\newcommand{\solver}[1]{\textsc{#1}}

\newcommand{\C}{\texttt{LC}}
\newcommand{\NC}{\texttt{NLC}}
\newcommand{\M}{\texttt{RF}\xspace}
\newcommand{\always}{\texttt{AlwaysLC}\xspace}
\newcommand{\never}{\texttt{NeverLC}\xspace}
\newcommand{\oracle}{\texttt{Oracle}}
\newcommand{\initbound}{c_\textrm{i}}
\newcommand{\dualbound}{c_\textrm{d}}
\newcommand{\primalbound}{c_\textrm{p}}
\newcolumntype{L}[1]{>{\raggedright\arraybackslash}p{#1}}
\newcolumntype{R}[1]{>{\raggedleft\arraybackslash}p{#1}}
\newcolumntype{C}[1]{>{\centering\arraybackslash}p{#1}}
\newcommand{\Met}{Met}
\newcommand{\miplib}[1]{\textsc{Miplib}~#1\xspace}


%% file: Features_Tab_twocol.tex
\begin{tabular}{|p{0.15cm}|ll|}  
\hline
& \textbf{Feature} & \textbf{Definition} \\
\hline
\hline
\multirow{26}{*}{\rotatebox[origin=c]{90}{\textnormal{\emph{\textbf{Static}}}}} & \multicolumn{2}{c|}{\emph{$\sim$ Matrix $\sim$}} \\
\hhline{~--}
& \texttt{Rows} & $\ln(m)$ \\
& \texttt{Columns} & $\ln(n)$ \\
& \texttt{NonZeros} & ratio of non-zeros, over $m\times n$ \\
& \texttt{Symmetries} & $1$ if any symmetry, $0$ otherwise \\
\hhline{~--}
& \multicolumn{2}{c|}{\emph{$\sim$ Variables $\sim$}} \\
\hhline{~--}
& \texttt{Binaries} & ratio of binary variables, over $n$ \\
& \texttt{Integers} & ratio of integer variables, over $n$ \\
\hhline{~--}
& \multicolumn{2}{c|}{\emph{$\sim$ Constraints $\sim$}} \\
\hhline{~--}
& \texttt{LessThan} & \multirow{16}{*}{ratio of constraints per constraint type, over $m$} \\
& \texttt{GreaterThan} & \\
& \texttt{Equality} & \\
& \texttt{SetPartitioning} & \\
& \texttt{SetPacking} & \\
& \texttt{SetCovering} &\\
& \texttt{Cardinality} & \\
& \texttt{KnapsackEquality} &\\
& \texttt{Knapsack} & \\
& \texttt{KnapsackInteger} & \\
& \texttt{BinaryPacking} & \\
& \texttt{VariableLowerBound} & \\
& \texttt{VariableUpperBound} & \\
& \texttt{MixedBinary} & \\
& \texttt{MixedInteger} & \\
& \texttt{Continuous} & \\
\hline
\multirow{12}{*}{\rotatebox[origin=c]{90}{\emph{\textbf{Dynamic}}}}  & \multicolumn{2}{c|}{\emph{$\sim$ Scaling $\sim$}} \\
\hhline{~--}
& \texttt{Coefficient\_Oom} & $\ln(\max A^{\prime} / \min A^{\prime})$ \\
& \texttt{RightHandSide\_Oom} & $\ln(\max b^{\prime} / \min b^{\prime})$ \\
& \texttt{Objective\_Oom} & $\ln(\max c^{\prime} / \min c^{\prime})$ \\
\hhline{~--}
& \multicolumn{2}{c|}{\emph{$\sim$ Presolving $\sim$}} \\
\hhline{~--}
& \texttt{PresolRows} & $\ln(\tilde{m})$ \\
& \texttt{PresolColumns} & $\ln(\tilde{n})$ \\
& \texttt{PresolIntegers} & ratio of presolved integer variables, over $n$ \\
\hhline{~--}
& \multicolumn{2}{c|}{\emph{$\sim$ Global Cutting $\sim$}} \\
\hhline{~--}
& \texttt{DualInitial\_Gap} & $\lvert\dualbound - \initbound\rvert/\max(\lvert \dualbound \rvert, \lvert \initbound \rvert, \lvert \dualbound-\initbound\rvert)$ \\ 
& \texttt{PrimalDual\_Gap} & $\lvert\primalbound - \dualbound\rvert/\max(\lvert\primalbound\rvert,\lvert \dualbound\rvert, \lvert\primalbound-\dualbound\rvert)$ \\
& \texttt{PrimalInitial\_Gap} & $\lvert\primalbound - \initbound\rvert/\max(\lvert\primalbound\rvert, \lvert\initbound\rvert, \lvert\primalbound-\initbound\rvert)$ \\
& \texttt{Gap\_Closed} & $1 - \frac{\texttt{PrimalDual\_Gap}}{\texttt{DualInitial\_Gap}}$ \\
\hline
\end{tabular}

%% file: Regression_Tab.tex
\begin{tabular}{|p{0.01cm}|L{0.8cm}R{1.15cm}R{1.15cm}R{1.15cm}R{0.7cm}R{1.15cm}R{0.7cm}R{0.95cm}|}
\hline
& \textbf{Metric} & $\M$ & \always & \never & \textbf{Imp.} ($\%$) & \oracle & \textbf{Pot.} ($\%$) & \textbf{Achiev.} ($\%$) \\
\hline
\multirow{3}{*}{\rotatebox[origin=c]{90}{\emph{Train}}}
& $\mathit{Time}$ & 427.91 & 467.92 & 642.56 & 8.55 & 414.75 & 11.36 & 75.25 \\  
& $\mathit{PDI}$ & 79.50 & 84.39 & 99.62 & 5.79 & 76.21 & 9.69 & 59.79 \\ 
& $\mathit{Nodes}$ & 22363.22 & 22693.97 & 43554.10 & 1.46 & 20363.27 & 10.27 & 14.19\\
\hline
\multirow{3}{*}{\rotatebox[origin=c]{90}{\emph{Test}}} 
& $\mathit{Time}$ & 420.09 & 434.61 & 593.21 & 3.34 & 384.72 & 11.48 & 29.09 \\  
& $\mathit{PDI}$ & 78.13 & 79.26 & 93.66 & 1.43 & 72.03 & 9.12 & 15.66 \\ 
& $\mathit{Nodes}$ & 21243.26 & 19903.08 & 37812.65 & -6.73 & 17919.92 & 9.96 & -67.58\\
\hline
\end{tabular}

%% file: Bracket_Tab.tex
\begin{tabular}{lrrrrrr}
\hhline{~~----}
\multicolumn{2}{c}{} & \multicolumn{2}{c}{\M} & \multicolumn{2}{c}{\always} & \multicolumn{1}{c}{} \\
\hline
\textbf{Bracket} & \textbf{Instances} & \textbf{Solved} & \textbf{Time} & \textbf{Solved} & \textbf{Time} & \textbf{Imp.} ($\%$) \\
\hline
All & 1155 & 872 & 426.59 & 854 & 462.16 &  7.70 \\
$[10, 7200]$ & 753 & 751 & 257.15 & 733 & 290.20 & 11.39 \\
$[100, 7200]$ & 470 & 468 & 750.73 & 450 & 874.37 & 14.14 \\
$[1000, 7200]$ & 226 & 224 & 2198.42 & 206 & 2703.86 & 18.69 \\
Affected & 310 & 308 & 132.45 & 290 & 180.63 & 26.67 \\
\hline
\end{tabular}

%% file: Bracket_Tab_Xpress.tex
\begin{tabular}{lrrrrrr}
\hhline{~~----}
\multicolumn{2}{c}{} & \multicolumn{2}{c}{\M} & \multicolumn{2}{c}{without \M} & \multicolumn{1}{c}{} \\
\hline
\textbf{Bracket} & \textbf{Instances} & \textbf{Solved} & \textbf{Time} & \textbf{Solved} & \textbf{Time} & \textbf{Imp.} ($\%$) \\
\hline
All             & 5530 & 5328  &  102.56 & 5319 &  104.04 & 1.44 \\
$[10, 7200]$    & 4434 & 4431  &  131.01 & 4422 &  133.22 & 1.69 \\
$[100, 7200]$   & 2226 & 2223  &  463.08 & 2214 &  474.39 & 2.44 \\
$[1000, 7200]$  &  574 &  571  & 2093.23 &  562 & 2203.17 & 5.25 \\                                  
Affected        & 488  & 477   & 116.79  & 468  & 137.80  & 17.99 \\
\hline
\end{tabular}